\documentclass[a4paper,12pt]{amsart}

\usepackage{amsmath}
\usepackage{amssymb}
\usepackage{mathrsfs}
\usepackage{enumerate}
\usepackage{ifthen}
\usepackage{graphicx}
\usepackage{caption}
\usepackage{subcaption}
\usepackage{float}
\usepackage[T1]{fontenc} %skandit

\nonstopmode \numberwithin{equation}{section}
\newtheorem{thm}{Theorem}[section]
\newtheorem{cor}{Corollary}[section]
\newtheorem{lem}{Lemma}[section]

\theoremstyle{definition}

\newtheorem{example}{Example}[section]

\newtheorem{rem}{Remark}[section]

\newcounter{minutes}
\divide\time by 60
\newcounter{hours}
\multiply\time by 60
\addtocounter{minutes}{-\time}
\newcounter {own}
\def\theown {\thesection       .\arabic{own}}

{\qed\bigskip}

\newcounter{alphabet}

\begin{document}

\title{On Harmonic Univalent Spirallike Mappings}

 \author{Md Firoz Ali*}
\address{Md Firoz Ali,
Department of Mathematics,
National Institute of Technology Durgapur,
713209, West Bengal, India.}
\email{firoz.ali@maths.nitdgp.ac.in}
\email{ali.firoz89@gmail.com}

\author{Sushil Pandit}
\address{Sushil Pandit,
Department of Mathema
tics,
National Institute of Technology Durgapur,
713209, West Bengal, India.}
\email{sushilpandit15594@gmail.com}


\def\thefootnote{}
\footnotetext{{\tiny File:~\jobname.tex,
printed: \number\year-\number\month-\number\day,
          \thehours.\ifnum\theminutes<10{0}\fi\theminutes }\\
*Corresponding author:~ Md Firoz Ali,~email:~ali.firoz89@gmail.com / firoz.ali@maths.nitdgp.ac.in
} \makeatletter\def\thefootnote{\@arabic\c@footnote}\makeatother

\begin{abstract}
In this article, we provide some necessary and sufficient coefficients conditions for a harmonic mapping to be hereditarily spirallike. Also, we give growth estimate for certain harmonic hereditarily spirallike mappings. Moreover, we connect the concept of harmonic  hereditarily spirallike mapping to analytic spirallike mapping and provide some examples in support of our results.
\end{abstract}

\thanks{}

\maketitle
\pagestyle{myheadings}
\markboth{Md Firoz Ali, Sushil Pandit}{On Harmonic Univalent Spirallike Mappings}

\section{Introduction}
\noindent
Let $\mathcal{A}$ denote the class of all analytic functions $h$ in the unit disk $\mathbb{D}=\{z\in\mathbb{C}:|z|<1\}$ of the form
\begin{align}\label{p3-001}
h(z)=z+\sum\limits_{n=2}^\infty a_nz^n,
\end{align}
and $\mathcal{S}$ be the class of univalent functions in $\mathcal{A}.$  An analytic function $h\in\mathcal{S}$ is called starlike if it maps the unit disk $\mathbb{D}$ onto a starlike domain with respect to origin i.e., if it satisfies ${\rm Re\,}\left(zh'(z)/h(z)\right)>0$ for $z\in\mathbb{D}\setminus\{0\}.$ Similarly, an analytic function $h\in\mathcal{S}$ is called convex if it maps the unit disk $\mathbb{D}$ onto a convex domain i.e., if it satisfies ${\rm Re\,}\left(1+zh''(z)/h'(z)\right)>0$ for $z\in\mathbb{D}\setminus\{0\}.$  A domain $\Omega$ containing the origin is called $\lambda-$ spirallike if for each point $w_0$ in $\Omega,$ the arc of $\lambda-$ spiral from origin to $w_0$ lies entirely in $\Omega.$ For $\lambda\in(-\pi/2,\pi/2),$ a function $h\in\mathcal{S}$ is called $\lambda-$ spirallike if $h(\mathbb{D})$ is a $\lambda-$ spirallike domain i.e., if ${\rm Re\,}\left(e^{-i\lambda}zh'(z)/h(z)\right)>0$ for $z\in\mathbb{D}\setminus\{0\}.$ For more information about these classes, we refer to \cite{Duren-1983}.\\

The differential operator $D$ for a continuously differentiable function $f\in C^1(\mathbb{D})$ is defined  by
$$Df(z)=zf_z(z)-\overline{z}f_{\overline{z}}(z),$$
where $f_z=(f_x-if_y)/2$ and $f_{\overline{z}}=(f_x+if_y)/2.$ Here $f_x$ and $f_y$ are the partial derivatives of $f$ with respect to $x={\rm Re\,}z$ and $y={\rm Im\,}z,$ respectively. A continuous twice differentiable complex valued function $f=u+iv$ is called harmonic in a complex domain $\Omega$ if both $u$ and $v$ are real harmonic in $\Omega.$ In any simply connected domain $\Omega,$ every harmonic mapping $f$ can be represented as $f=h+\overline{g}$ where $h$ and $g$ are analytic in $\Omega.$ This representation is known as canonical representation and $h$ is called analytic part whereas $g$ is called co-analytic part of $f.$ A harmonic mapping $f=h+\overline{g}$ is sense-preserving if the Jacobian $J_f=|h'|^2-|g'|^2$ is positive and sense-reversing if $J_f$ is negative. Lewy \cite{Lewy-1936} showed that $f$ is locally univalent if $J_f$ is non vanishing. Let $\mathcal{H}$ be class of all sense-preserving harmonic functions $f=h+\overline{g}$ in $\mathbb{D}$ with
\begin{align}\label{p3-010}
h(z)=z+\sum\limits_{n=2}^\infty a_nz^n,~g(z)=\sum\limits_{n=1}^\infty b_nz^n
\end{align}
and $\mathcal{H}'$ be class of harmonic functions $f=h+\overline{g}$ in $\mathcal{H}$ with the representation
\begin{align}\label{p3-020}
 h(z)=z-\sum\limits_{n=2}^\infty |a_n|z^n, ~g(z)=\sum\limits_{n=1}^\infty |b_n|z^n.
\end{align}
Let $\mathcal{S}_H$ denote the family of functions in $\mathcal{H}$ that are univalent. In 1984, Clunie and Sheil-Small \cite{Clunie-Small-1984} studied several geometric properties of functions in $\mathcal{S}_H$ alongwith functions in its  subclasses of starlike functions, convex functions, close-to-convex functions,  etc.\\

It is well known that the properties convexity and starlikeness of an analytic and univalent mapping are hereditary i.e., if an analytic and univalent function $f$ is convex or starlike then $f(\mathbb{D}_r)$ is a convex or starlike domain, respectively for each $0<r<1,$ where $\mathbb{D}_r=\{z\in\mathbb{C}:|z|<r\}.$ But the properties convexity and starlikeness are not hereditary for a harmonic univalent map $f$. For example, the half-plane harmonic mapping $l(z)=h(z)+\overline{g(z)}$ with
\begin{align*}
h(z)=\frac{z-\frac{1}{2}z^2}{(1-z)^2},\quad g(z)=\frac{-\frac{1}{2}z^2}{(1-z)^2}
\end{align*}
maps $\mathbb{D}$ univalently onto ${\rm Re\,}w >-1/2$, whereas $l(\mathbb{D}_r)$ is not convex for $\sqrt{2}-1<r<1$ (see \cite{Ma-Ponnusamy-Sugawa-2022}). On the other hand, the harmonic Koebe function $k(z)=h(z)+\overline{g(z)}$ with
\begin{align*}
h(z)=\frac{z-\frac{1}{2}z^2+\frac{1}{6}z^3}{(1-z)^3},\quad g(z)=\frac{\frac{1}{2}z^2+\frac{1}{6}z^3}{(1-z)^3}
\end{align*}
map $\mathbb{D}$ onto a slit domain $\mathbb{C}\setminus(-\infty,-1/6)$ in one-to-one manner (see \cite{Clunie-Small-1984}) but, recently, Ma et al. \cite{Ma-Ponnusamy-Sugawa-2022} showed that $k(\mathbb{D}_r)$ is not starlike for $r=\sqrt{5}/3.$ Due to this phenomena, in 2004, Chuaqui et al. \cite{Chuaqui-Duren-Osgood-2004} introduced the concept of fully convexity and fully starlikeness for harmonic mappings. A harmonic mapping $f=h+\overline{g}\in\mathcal{H}$ is called fully convex if it maps every circle $C_r=\{z\in\mathbb{C}:|z|=r\},~0<r<1,$  in an one-to-one manner onto a convex curve. It is well known that, a sense-preserving harmonic mapping $f=h+\overline{g}\in\mathcal{H}$ is fully convex if ${\rm Re\,}\left(\frac{D^2f(z)}{Df(z)}\right)>0$ for every $z=re^{i\theta}\in\mathbb{D}\setminus\{0\}.$ Similarly, a harmonic mapping $f=h+\overline{g}\in\mathcal{H}$ is said to be fully starlike if it maps every circle $C_r,~0<r<1,$ in an one-to-one manner onto a curve that bounds a domain starlike with respect to the origin. It is well known that, a sense-preserving harmonic mapping $f=h+\overline{g}\in\mathcal{H}$ with $f(0)=0$ if and only if $z=0$ is fully starlike if ${\rm Re\,}\left(\frac{Df(z)}{f(z)}\right)>0$ for every $z=re^{i\theta}\in\mathbb{D}\setminus\{0\}.$ In the same paper \cite{Chuaqui-Duren-Osgood-2004}, authors showed that  a harmonic  fully convex function $f=h+\overline{g}\in\mathcal{H}$ is also univalent but a harmonic fully starlike function $f=h+\overline{g}\in\mathcal{H}$ may not be even locally univalent. As in \cite{Ma-Ponnusamy-Sugawa-2022}, we call a harmonic fully starlike and univalent mapping $f=h+\overline{g}\in\mathcal{H}$ as harmonic hereditarily starlike function. Let $\mathcal{ST}_H$ denotes the class of harmonic  mappings $f\in\mathcal{H}$ which are hereditarily starlike and $\mathcal{ST}_{H'}$ be the class of harmonic mappings $f\in\mathcal{ST}_H$ of the form \eqref{p3-020}. Different geometric aspects of starlikeness and hereditarily starlikeness for harmonic mappings in the class $\mathcal{H}$ and $\mathcal{H}'$ have been studied in \cite{Chuaqui-Duren-Osgood-2004, Jahangiri-1999, Mocanu-1980, Nagpal-Ravichandran-2013, Ponnusamy-Kaliraj-2014}.\\

Clunie and Sheil-Small \cite{Clunie-Small-1984} successfully extended the theory of starlikeness, convexity and close-to-convexity from analytic functions to harmonic functions. So it is natural to think whether spirallkeness can be extended from analytic functions to harmonic functions. Affirmatively, Ma et al. \cite{Ma-Ponnusamy-Sugawa-2022} introduced the concept of hereditarily spirallikeness for harmonic mappings in the unit disk $\mathbb{D}$ and investigated its geometric properties. For a real number $\lambda$ with $|\lambda|<\frac{\pi}{2}$, a harmonic function $f=h+\overline{g}$ in $\mathcal{H}$ is called hereditarily  $\lambda-$spirallike if $f$ is univalent on $\mathbb{D}$ and $f(\mathbb{D}_r)$ is a $\lambda-$spirallike domain for each $r<1.$ The following analytic characterization of hereditarily spirallike functions has been obtained in \cite{Ma-Ponnusamy-Sugawa-2022}.

\begin{lem}\cite{Ma-Ponnusamy-Sugawa-2022} \label{p3-030}
Let $\lambda\in (-\pi/2,\pi/2)$. Suppose that a function $f\in C^1(\mathbb{D})$ satisfies the condition that $f(z)=0$ if and only if $z=0,$ and that $J_f=|f_z|^2-|f_{\overline{z}}|^2>0$ on $\mathbb{D}.$ Then $f$ is one-one in $\mathbb{D}$ and $f(\mathbb{D}_r)$ is $\lambda-$spirallike for each $0<r<1$ if and only if
\begin{align*}
{\rm Re\,}\left(e^{-i\lambda}\frac{Df(z)}{f(z)}\right)>0,\quad z\in\mathbb{D}\setminus\{0\}.
\end{align*}
\end{lem}
Let $\mathcal{SP}_H(\lambda)$ and $\mathcal{SP}_{H'}(\lambda)$ denote the classes of harmonic hereditarily  $\lambda-$spirallike functions in $\mathcal{H}$ and $\mathcal{H}',$ respectively. \\

%There are many other articles on univalent harmonic mapping defined in the unit disk $\mathbb{D}=\{z\in\mathbb{C}: |z|<1\}.$  For convenient, we refer a few \cite{Chen-Gauthier-Hengartner-2000, Liu-Ponnusamy-2018, Silverman-1998}. The study of harmonic mappings in the unit disk $\mathbb{D}$ has been eye catching area in complex analysis.
%\noindent
In this article, we provide some necessary and sufficient conditions for a harmonic mapping to be hereditarily spirallike function. Also, we give growth estimate for certain harmonic hereditarily spirallike mappings. Moreover, we connect this concept of harmonic  hereditarily spirallike mapping with analytic spirallike function and  develop some examples in support of our results. In Section \ref{Background And Main Results}, we present our results and examples. In Section \ref{Proof of Main Results}, we provide detail proof of the results.
% Before that, we collect some notations which we use throughout the manuscript. For $0\leq \alpha<1,$ let $\mathcal{ST}_H(\alpha)$ denote the class of harmonic  mappings $f\in\mathcal{H}$ which are univalent and fully starlike of order $\alpha$ in $\mathbb{D}$ and $\mathcal{ST}_H'(\alpha)$ be class of mappings $f=h+\overline{g}$ in $\mathcal{ST}_H(\alpha)$ of the form \eqref{p3-020}. Throughout the paper we use $\mathcal{ST}_H(0)=\mathcal{ST}_H$ and $\mathcal{ST}_H'(0)=\mathcal{ST}_H'.$   The class of harmonic functions $f=h+\overline{g}\in\mathcal{H}$ which are univalent and hereditarily  $\lambda-$spirallike  will be denoted by $\mathcal{SP}_H(\lambda)$ and  let $\mathcal{SP}_{H'}(\lambda)$ denote the class of harmonic mappings $f=h+\overline{g}\in\mathcal{SP}_H(\lambda)$ of the form \eqref{p3-020}.
%\begin{align}\label{pp-55}
%h(z)=z-\sum\limits_{n=2}^\infty |a_n|z^n,~g(z)=\sum\limits_{n=1}^\infty |b_n|z^n
%\end{align}
%by

\section{Main Results}\label{Background And Main Results}
In 2004, Chuaqui et al. \cite{Chuaqui-Duren-Osgood-2004} obtained a necessary and sufficient coefficient condition for harmonic fully starlike and harmonic fully convex  mappings. In this context, we obtain a similar result for harmonic hereditarily spirallike  mappings.

\begin{thm}\label{p3-040}
Let $f=h+\overline{g}$ be a sense-preserving harmonic mapping on $\mathbb{D}$ such that $f(0)=0$ only for $z=0.$ Then $f\in\mathcal{SP}_H(\lambda)$ if and only if
\begin{align}\label{p3-050}
|h(z)|^2{\rm Re\,}\left(e^{-i\lambda}\frac{zh'(z)}{h(z)}\right)
&>|g(z)|^2{\rm Re\,}\left(e^{i\lambda}\frac{zg'(z)}{g(z)}\right)\\
&\qquad\quad +{\rm Re\,}\left(ze^{i\lambda}(h(z)g'(z)-e^{-2i\lambda}g(z)h'(z))\right)\nonumber.
\end{align}
\end{thm}
Theorem \ref{p3-040} can be used to identify whether a harmonic mapping is hereditarily spirallike or not. In this context, we present the following example.
\begin{example}\label{p3-060}
For $\alpha\in\mathbb{D},$ the harmonic mapping $f_1(z)=z+\alpha\overline{z}$ is not hereditarily $\pi/4-$spirallike. Because for $\alpha=-1/2$ and $z=1/2(1+i)$ the inequality \eqref{p3-050} becomes $3/4>1.$
\end{example}

Here we note that, for $\alpha\in\mathbb{D}$ and an analytic univalent mapping $h\in\mathcal{S},$ the harmonic mapping $f=h+\alpha\overline{h}$ is sense-preserving fully starlike in $\mathbb{D}$ if and only if $h$ is starlike in $\mathbb{D},$ because
$${\rm Re\,}\left(\frac{Df(z)}{f(z)}\right)=\frac{(1-|\alpha|^2)|h(z)|^2}{|h(z)+\alpha\overline{h(z)}|^2}{\rm Re\,}\left(\frac{zh'(z)}{h(z)}\right).$$
But the Example \ref{p3-060} shows that this does not happen in the case of spirallike mappings, and so the class $\mathcal{SP}_H(\lambda)$  is not affine invariant.\\
% by showing that $h(z)+\alpha\overline{h(z)}$ is not harmonic hereditarily $\lambda-$spirallike for some analytic $\lambda-$spirallike function $h(z).$ For example, we consider $h(z)=z$ which is $\lambda-$spirallike for $\lambda\in(-\pi/2,\pi/2)$ but, later in the next section, we will show that $f(z)=z+\alpha\overline{z}$ is not hereditarily $\pi/4-$spirallike for $\alpha=-1/2.$

%Firstly, we note down a recent result by Xiu-Shuang Ma et.al \cite{Ma-Ponnusamy-Sugawa-2022} which is helpful to handle the proof.
%[\textbf{Proof of Theorem \ref{p3-040}}]
%\begin{proof}

%If $f\in\mathcal{SP}_H(\lambda)$ then ${\rm Re\,}\left(e^{-i\lambda}\frac{Df(z)}{f(z)}\right)>0$ for $z\in\mathbb{D}$. using the fact ${\rm Re\,}w={\rm Re\,}(\overline{w})$ for $w\in\mathbb{C}$ in the last inequality, one can easily get \eqref{p3-050}. \\
%Conversely, a simple calculation leads the condition \eqref{p3-050} to ${\rm Re\,}\left(e^{-i\lambda}\frac{Df(z)}{f(z)}\right)>0.$ Furthermore, $f$ is sense-preserving and $f(z)=0$ only for $z=0.$ Therefore Lemma \ref{p3-030} implies that the function $f$ is hereditarily $\lambda-$spirallike.
%\end{proof}

In 1998, Silverman \cite{Silverman-1998} considered the classes $\mathcal{H}$ and $\mathcal{H}'$ of harmonic mappings and obtained sufficient conditions  for functions in $\mathcal{H}$ and $\mathcal{H}'$ to be a member of the class $\mathcal{ST}_H$ and $\mathcal{ST}_{H'},$ respectively, in term of coefficients.  Furthermore, the author proved that the obtained condition is also necessary for the the later case. For convenient of the reader, we state the result here.
\begin{lem}\label{p3-070}\cite{Silverman-1998}
Let $f=h+\overline{g}\in\mathcal{H}$ be of the form \eqref{p3-010} such that
\begin{align}\label{p3-080}
\sum\limits_{n=1}^\infty n\left(|a_n|+|b_n|\right)\leq 2.
\end{align}
Then $f\in\mathcal{ST}_H.$ Moreover, the condition \eqref{p3-080} is necessary and sufficient for a harmonic mapping  $f=h+\overline{g}\in\mathcal{H}'$ of the form \eqref{p3-020} to be a member of the class $\mathcal{ST}_H'.$
\end{lem}
%\begin{lem}\label{pp-59}
%Let $f=h+\overline{g}$ be given by \eqref{pp-55}. Then $f\in\mathcal{ST}_H(\alpha)$ if and only if
%\begin{align*}
%\sum\limits_{n=1}^\infty\left(\frac{n-\alpha}{1-\alpha}|a_n|+\frac{n+\alpha}{1-\alpha}|b_n|\right)\leq 2.
%\end{align*}
%\end{lem}
 We present a few similar results for harmonic hereditarily spirallike mappings in the unit disk. Before that we note down a few notations which we will use throughout the paper. For $\lambda\in(-\pi/2,\pi/2)$ and $n\geq 1,$ let
 \begin{align}\label{p3-090}
 A_n=\left|1+ne^{-i\lambda}\right|+\left|1-ne^{-i\lambda}\right|\quad\text{and}\quad B=\left|1+e^{-i\lambda}\right|-\left|1-e^{-i\lambda}\right|.
\end{align}

 % =(1+n^2+2n\cos\lambda)^{1/2}+(1+n^2-2n\cos\lambda)^{1/2}\\
%=\sqrt{2}[(1+\cos\lambda)^{1/2}-(1-\cos\lambda)^{1/2}]
 \begin{thm}\label{p3-100}
Let $\lambda\in(-\pi/2,\pi/2)$ and $f=h+\overline{g}$ be of the form \eqref{p3-010} such that
\begin{align}\label{p3-110}
\sum\limits_{n=2}^\infty\frac{A_n}{B}|a_n|+\sum\limits_{n=1}^\infty\frac{A_n}{B}|b_n|\leq 1
\end{align}
where $A_n$ and $B$ are given by \eqref{p3-090}. Then $f$ is harmonic hereditarily $\lambda-$spirallike function.
\end{thm}

For $\lambda=0,$ Theorem \ref{p3-100} reduces to Lemma \ref{p3-070}, which is expected because $\mathcal{SP}_H(0)=\mathcal{ST}_H.$ \\

There are harmonic hereditarily $\lambda-$spirallike mappings for which the value of the left hand side of \eqref{p3-110} is $1.$ For example, the harmonic hereditarily $\lambda-$spirallike mappings
\begin{align*}
f(z)=h(z)+\overline{g(z)}=z+\sum\limits_{n=2}^\infty\frac{B}{A_n}x_nz^n+\sum\limits_{n=1}^\infty\frac{B}{A_n}\overline{y_nz^n}
\end{align*}
where $\sum\limits_{n=2}|x_n|+\sum\limits_{n=1}|y_n|=1$ have this  property.\\
%Here, we notice that Lemma \ref{p3-070} is a particular case of our result in Theorem \ref{p3-100}. We have highlighted this in the following corollary which can be easily verified.
%\begin{cor}
%Theorem \ref{p3-100} coincides the Lemma \ref{p3-070} for $\lambda=0.$
%\end{cor}
%%\begin{cor}
%%g
%%\end{cor}
%\begin{rem}
%From the inequality \eqref{p3-110}, it follows that modulus of the dilatation at origin is $|\omega(0)|=|b_1|\leq B/A_1.$
%\end{rem}

\begin{figure}[H]
\subfloat[$\alpha=0.95$ and $\lambda=\pi/3$\label{sp3}]{%
\includegraphics[width=0.43\textwidth, height=7cm]
{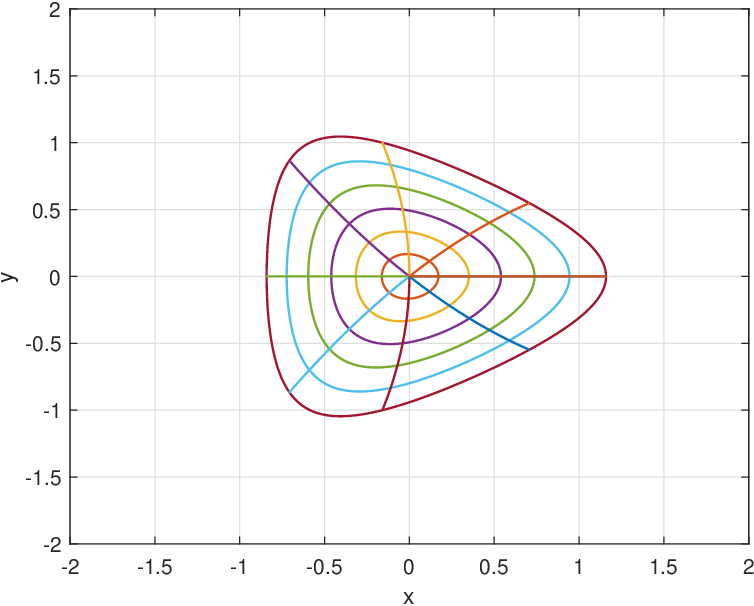}}
\hspace{.2cm}
\subfloat[$\alpha=0.95$ and $\lambda=\pi/5$\label{sp5}]{%
\includegraphics[width=0.43\textwidth, height=7cm]
{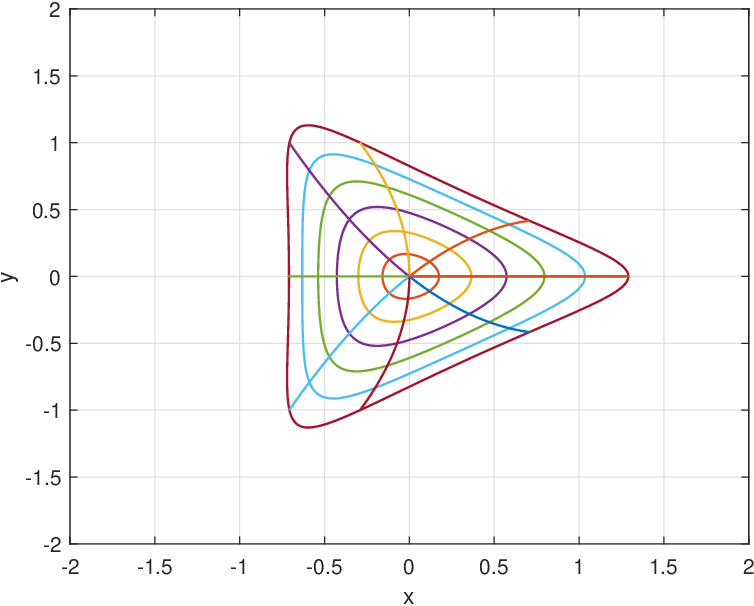}}
\\
\vspace{.2cm}
\subfloat[$\alpha=0.95$ and $\lambda=\pi/8$\label{sp8}]{%
\includegraphics[width=0.43\textwidth, height=7cm]
{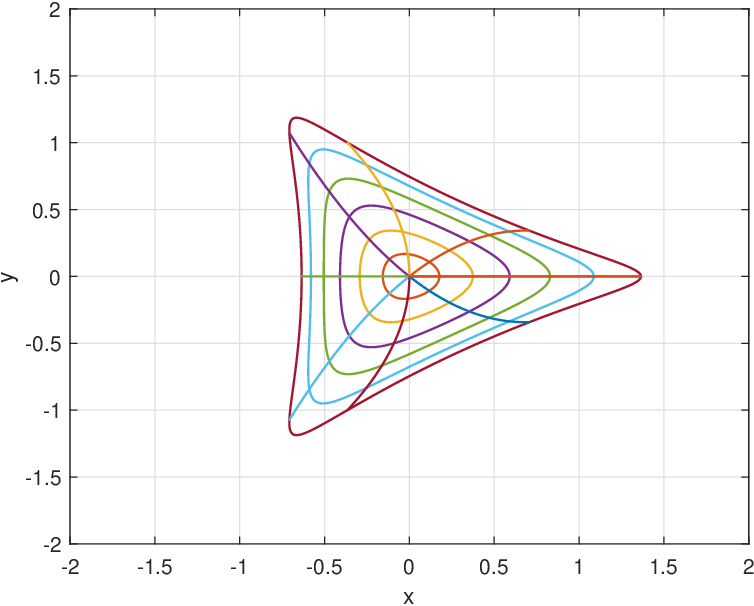}}
\hspace{.2cm}
\subfloat[$\alpha=0.95$ and $\lambda=0$\label{sp0}]{%
\includegraphics[width=0.43\textwidth, height=7cm]
{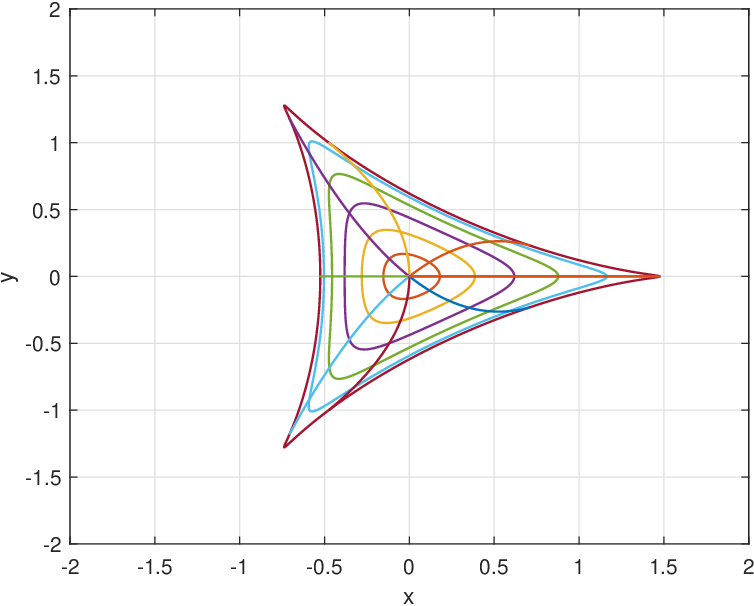}}
\caption{$f_2(\mathbb{D})$ for certain values of $\alpha$ and $\lambda.$}
\label{fig-1}
\end{figure}

The inequality \eqref{p3-110} is useful to construct harmonic univalent spirallike functions in the unit disk. Using the inequality \eqref{p3-110}, one can easily check that $f_2(z)=z+\alpha B/A_2\overline{z^2}$  and  $f_3(z)=z+B\alpha/A_1\overline{z}+B/A_3(1-|\alpha|)\overline{z^3}$ where $\alpha\in\mathbb{D}$ are sense-preserving harmonic hereditarily $\lambda-$spirallike mappings for every $\lambda\in(-\pi/2,\pi/2)$. The images of $\mathbb{D}$ under $f_2$ for certain values of $\alpha$ and $\lambda$ are shown in Figure \ref{fig-1}.

%  \begin{figure}[H]
%   \includegraphics[width=0.5cm]{sp3}\hfill
%     \caption{For $\alpha=.95$ and $\lambda=\pi/3$}
%   \includegraphics[width=0.5cm]{sp5}\hfill
%     \caption{For $\alpha=.95$ and $\lambda=\pi/5$}
%  \includegraphics[width=0.5cm]{sp12}
%  \caption{For $\alpha=.95$ and $\lambda=\pi/12$}
%\end{figure}
%\noindent
 The condition \eqref{p3-110} in the Theorem \ref{p3-100} is not a necessary condition. This can be seen by an example of a hereditarily $\lambda-$spirallike function which does not satisfy the condition \eqref{p3-110}. We consider the function $f_4(z)=z(1-z)^{i-1}=z+\sum\limits_{n=2}^\infty a_nz^n\in\mathcal{H}.$ It is known that $f_4$ is hereditarily $\pi/4-$spirallike function (see \cite{Robertson-1969}). Clearly,
 $$a_2=\frac{f_3''(0)}{2}=1-i.$$
But for $\lambda=\pi/4,$ we have
 $$\sum\limits_{n=2}^\infty\frac{A_n}{B}|a_n|+\sum\limits_{n=1}^\infty\frac{A_n}{B}|b_n|\geq \frac{A_2}{B}|a_2|=\sqrt{2}\frac{A_2}{B}>2.$$

 Our next result gives necessary and sufficient condition for harmonic mappings in $\mathcal{H}'$ to be a member of the class $\mathcal{SP}_{H'}(\lambda).$
\begin{thm} \label{p3-120}
Let $\lambda\in(-\pi/2,\pi/2)$ and $f=h+\overline{g}\in\mathcal{H}'$ be of the form \eqref{p3-020}. If $f$ satisfies \eqref{p3-110}
%\begin{align*}
%\sum\limits_{n=1}^\infty\frac{A_n}{B}\left(|a_n|+|b_n|\right)\leq 2.
%\end{align*}
then $f\in\mathcal{SP}_{H'}(\lambda).$ Conversely, every harmonic function $f\in\mathcal{SP}_{H'}(\lambda)$ satisfies the coefficient inequality
\begin{align}\label{p3-130}
\sum\limits_{n=2}^\infty\frac{B}{A_n}|a_n|+\sum\limits_{n=1}^\infty\frac{B}{A_n}|b_n|\leq 1
\end{align}
where $A_n$ and $B$ are given by \eqref{p3-090}.
\end{thm}

%Again, the condition \eqref{p3-110} in the Theorem \ref{p3-120} is not necessary. For example we consider a harmonic mapping $f_4(z)=z+\alpha\overline{z}\in\mathcal{H}'$ for $0\leq \alpha<1.$ A simple calculation shows that
%$${\rm Re\,}\left(e^{-i\lambda}\frac{Df(z)}{f(z)}\right)>0$$ is equivalent to
%\begin{align}\label{pp-175}
%(1-|\alpha|^2)|z|^2>2\tan{\lambda}{\rm Im\,}(\alpha z^2)~\text{for} \lambda\in[0,\pi/2).
%\end{align}
%By the fact ${\rm Im\,}(\alpha z^2)\leq|\alpha||z|^2,$ we say that \eqref{pp-175} holds if
%$$(1-|\alpha|^2)>|\alpha|\tan{\lambda}~\text{for} \lambda\in[0,\pi/2).$$
%In particular, for $\alpha$ for $\alpha=5/6\sqrt{2}<1/\sqrt{2}$ the function $f_4$ is hereditarily $\pi/4-$spirallike but the value of
%$$\sum\limits_{n=2}^\infty\frac{A_n}{B}|a_n|+\sum\limits_{n=1}^\infty\frac{A_n}{B}|b_n|\geq\frac{A_1|b_1|}{B}\thickapprox 1.42258898>1.$$

The proof of the converse part of the Theorem \ref{p3-120} provides the following result.
\begin{cor}\label{p3-140}
Let $\lambda\in(-\pi/2,\pi/2)$ and $f=h+\overline{g}\in\mathcal{SP}_{H'}(\lambda)$ be of the form \eqref{p3-020}. Then $\sum\limits_{n=1}^\infty(n|a_n|+n|b_n|)\leq 2,$ and $f\in\mathcal{ST}_H'.$
\end{cor}

By Theorem \ref{p3-100} and Theorem \ref{p3-120}, we have constructed harmonic hereditarily $\lambda-$spirallike mapping in infinite series form.
\begin{example}\label{p3-150}
Let $h_1(z)=z,~ h_n(z)=z+B/A_nz^n,~n\geq2$ and $g_n(z)=z+B/A_n\overline{z^n},~n\geq1$ and $X_n,Y_n\geq0$ are such that $\sum\limits_{n=1}^\infty(X_n+Y_n)=1.$ Then $f(z)=\sum\limits_{n=1}^\infty(X_nh_n(z)+Y_ng_n(z))$ is hereditarily $\lambda-$spirallike.\\

Clearly, $f$ can be written as
\begin{align*}
f(z)&=\sum\limits_{n=1}^\infty(X_nh_n(z)+Y_ng_n(z))\\
    &=\sum\limits_{n=1}^\infty(X_n+Y_n)z+\sum\limits_{n=2}^\infty\left(\frac{X_nB}{A_n}z^n\right)+\sum\limits_{n=1}^\infty\left(\frac{Y_nB}{A_n}\overline{z^n}\right)\\
    &=z+\sum\limits_{n=2}^\infty\left(\frac{X_nB}{A_n}z^n\right)+\sum\limits_{n=1}^\infty\left(\frac{Y_nB}{A_n}\overline{z^n}\right).
\end{align*}
Thus the function $f$ is of the form \eqref{p3-010} and satisfies
\begin{align*}
\sum\limits_{n=2}^\infty\frac{A_n}{B}\left(\frac{X_nB}{A_n}\right)+\sum\limits_{n=1}^\infty\frac{A_n}{B}\left(\frac{Y_nB}{A_n}\right)=\sum\limits_{n=2}^\infty X_n+\sum\limits_{n=1}^\infty Y_n=1-X_1\leq1.
\end{align*}
Hence by the Theorem \ref{p3-100}, $f$ is hereditarily $\lambda-$spirallike.
\end{example}
\begin{rem}
If we take $h_n(z)=z-B/A_nz^n,~n\geq2$ in the Example \ref{p3-150}, then  the infinite series $\sum\limits_{n=1}^\infty(X_nh_n(z)+Y_ng_n(z))$ ultimately belong  to the class $\mathcal{SP}_{H'}(\lambda).$
\end{rem}
%Our next result gives an infinite series representation for functions in the class $\mathcal{SP}_{H'}(\lambda).$
\begin{example}\label{p3-160}
Let $f=h+\overline{g}\in\mathcal{SP}_{H'}(\lambda).$ Then $f$ can be written as $f(z)=\sum\limits_{n=1}^\infty(X_nh_n(z)+Y_ng_n(z)),$ where
$h_1(z)=z,~ h_n(z)=z-B/A_nz^n,~n\geq2$ and $g_n(z)=z+B/A_n\overline{z^n},~n\geq1$ with some $X_n,Y_n\geq0$ are such that $\sum\limits_{n=1}^\infty(X_n+Y_n)=1.$

Since $f=h+\overline{g}\in\mathcal{SP}_{H'}(\lambda)$ is of the form \eqref{p3-020}, by Theorem \ref{p3-120} the condition \eqref{p3-130} holds. Let
\begin{align*}
h_1(z)=z,~h_n(z)=z-\frac{A_n}{B}z^n,~n\geq2\quad\text{and}\quad g_n(z)=z+\frac{A_n}{B}\overline{z^n},~n\geq1
\end{align*}
and set
\begin{align*}
X_n=\frac{B|a_n|}{A_n},~n\geq2\quad\text{and}\quad Y_n=\frac{B|b_n|}{A_n},~n\geq1
\end{align*}
with $X_1=1-\sum\limits_{n=2}^\infty X_n-\sum\limits_{n=1}^\infty Y_n.$ Then from  the condition \eqref{p3-130}, it follows that $0\leq X_n,~Y_n\leq1$ for $n\geq1.$
%Moreover,
% $$\sum\limits_{n=1}^\infty (X_n+Y_n)=1.$$ \\
Therefore,
\begin{align*}
\sum\limits_{n=1}^\infty(X_nh_n(z)+Y_ng_n(z))&=\sum\limits_{n=1}^\infty (X_n+Y_n)z-\sum\limits_{n=2}^\infty X_n\frac{A_n}{B}z^n+\sum\limits_{n=1}^\infty Y_n\frac{A_n}{B}\overline{z^n}\\
&=z-\sum\limits_{n=2}^\infty |a_n|z^n+\sum\limits_{n=1}^\infty |b_n|\overline{z^n}=f(z).
\end{align*}
\end{example}

The growth estimate for sense-preserving  harmonic $\lambda-$spirallike functions is given below.
\begin{thm}\label{p3-170}
Let $f=h+\overline{g}\in\mathcal{H}$ be a harmonic mapping of the form \eqref{p3-010} which satisfies the condition \eqref{p3-110}. Then the sharp inequality
\begin{align*}
\left(1-\frac{B}{A_1}\right)r \leq |f(z)| \leq \left(1+\frac{B}{A_1}\right)r \quad\text{for}~|z|=r
\end{align*}
\end{thm}
holds and
%\begin{cor}
%Let $f=h+\overline{g}\in\mathcal{H}$ satisfies \eqref{p3-110}. Then
\begin{align*}
\left\{w\in\mathbb{C}:|w|<1-\frac{B}{A_1}\right\}\subset f(\mathbb{D}).
\end{align*}
%\end{cor}
%%%Many geometric properties of univalent analytic function are known so it is quite interesting to connect harmonic function theory to univalent function theory as it provide more exposures.

There is an interesting relation between analytic starlike and spirallike functions which was obtained by Ba\c{s}g\"{o}ze and Keogh \cite{Basgoze-Keogh-1970}. For convenient, we state the result below.
\begin{lem}\label{p3-180}\cite{Basgoze-Keogh-1970}
Let $\lambda\in(-\pi/2,\pi/2).$ For each analytic  $\lambda-$spirallike function $h\in\mathcal{A}$ there exist an unique analytic starlike function $g\in\mathcal{A}$ such that
\begin{align}\label{p3-190}
\frac{h(z)}{z}=\left(\frac{g(z)}{z}\right)^{e^{-i\lambda}\cos\lambda}.
\end{align}
\end{lem}
We extend the Lemma \ref{p3-180} for functions in the class $C^1(\mathbb{D}).$
\begin{thm}\label{p3-200}
Let $\lambda$ be a real number with $|\lambda|<\pi/2.$ For a  hereditarily $\lambda-$spirallike function $h\in C^1(\mathbb{D})$ satisfying $h(z)=0$ only for $z=0$ and $J_h=|h_z|^2-|h_{\overline{z}}|^2>0$ for $z\in\mathbb{D}$ there exists a fully starlike function $g\in C^1(\mathbb{D})$ satisfying $g(z)=0$ only for $z=0$ and $J_g=|g_z|^2-|g_{\overline{z}}|^2>0$ for $z\in\mathbb{D}$ such that relation \eqref{p3-190} is satisfied.
%\begin{align}
%\frac{h(z)}{z}=\left(\frac{g(z)}{z}\right)^{e^{-i\lambda}\cos\lambda}
%\end{align}
\end{thm}
It is important to note that, the relation \eqref{p3-190} hold for analytic functions, but not for complex-valued harmonic functions as the exponential power of harmonic function may not be harmonic. For example, for the harmonic function $g(z)=z+\alpha\overline{z},~\alpha\in\mathbb{D},$  the corresponding mapping $h$ define by the relation \eqref{p3-190}, is not harmonic. Now, we give a relation between harmonic spirallike and starlike functions in some different way.

\begin{thm}\label{p3-210}
Let $F(z)=z-\sum\limits_{n=2}^\infty |a_n|z^n+\overline{\sum\limits_{n=1}^\infty |b_n|z^n}$ be harmonic hereditarily starlike mapping. Then $f(z)=z+\sum\limits_{n=2}^\infty d_na_nz^n+\overline{\sum\limits_{n=1}^\infty d_nb_nz^n},$ where $\{d_n\}$ is a sequence such that $|d_n|\leq nB/A_n$ for $n\geq 1,$ is harmonic hereditarily $\lambda-$spirallike mapping. Conversely, if $f(z)=z-\sum\limits_{n=2}^\infty |a_n|z^n+\overline{\sum\limits_{n=1}^\infty |b_n|z^n}$ is harmonic hereditarily $\lambda-$spirallike mapping then $F(z)=z+\sum\limits_{n=2}^\infty a_nz^n+\overline{\sum\limits_{n=1}^\infty b_nz^n}$ is hereditarily starlike mapping.
\end{thm}

\begin{figure}[H]
\subfloat[$\alpha=0.2$ and $\lambda=\pi/3$]{%
\includegraphics[width=0.4\textwidth, height=5cm]
{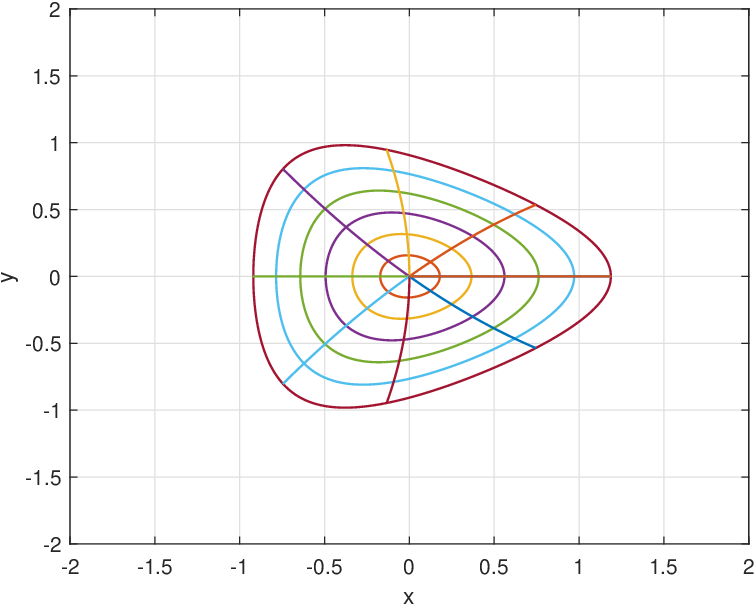}}
\hspace{.2cm}
\subfloat[$\alpha=0.5$ and $\lambda=\pi/3$]{%
\includegraphics[width=0.4\textwidth, height=5cm]
{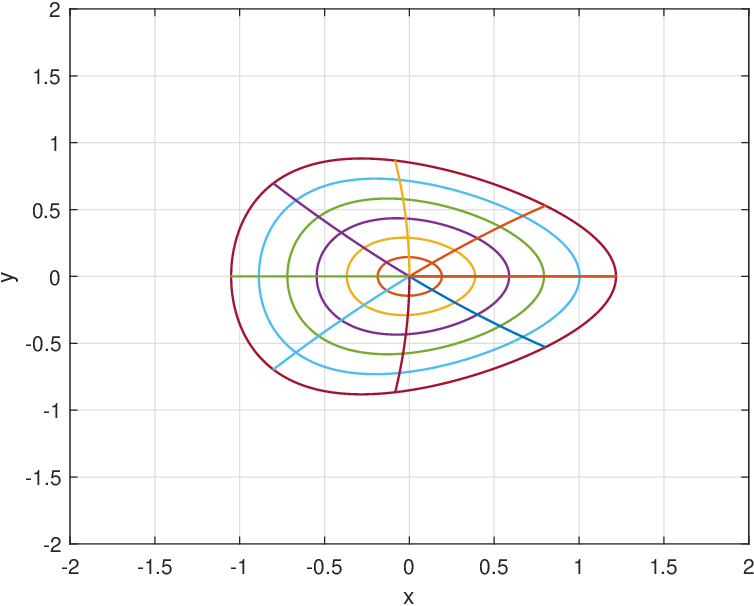}}
\\
\vspace{.2cm}
\subfloat[$\alpha=0.7$ and $\lambda=\pi/3$]{%
\includegraphics[width=0.4\textwidth, height=5cm]
{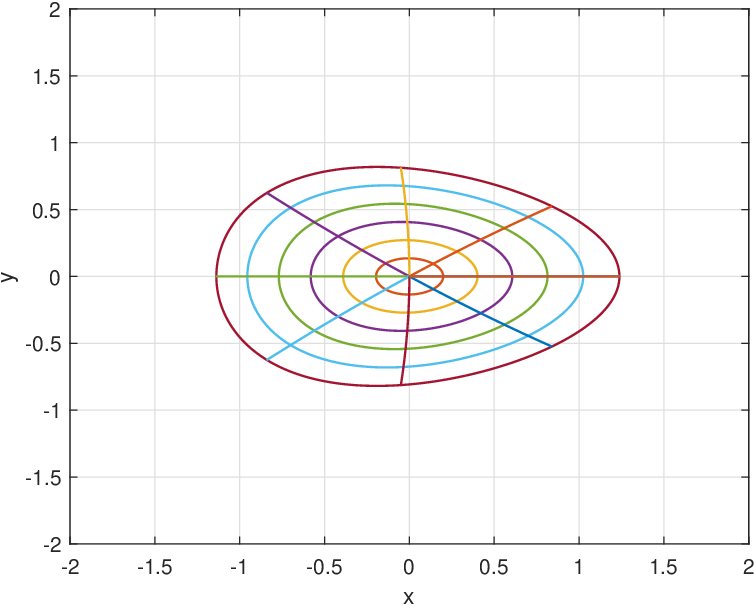}}
\hspace{.2cm}
\subfloat[$\alpha=0.9$ and $\lambda=\pi/3$]{%
\includegraphics[width=0.4\textwidth, height=5cm]
{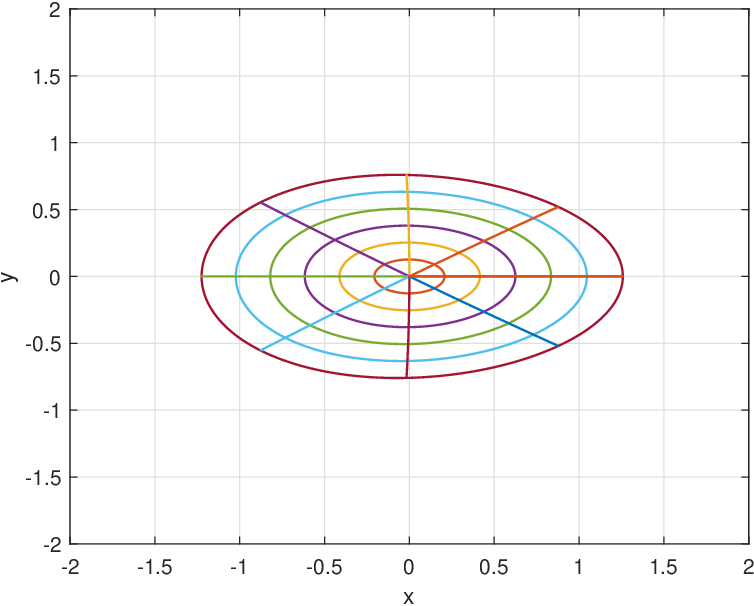}}
\caption{$f_5(\mathbb{D})$ for certain values of $\lambda$ and $\alpha.$}
\label{fig-3}
\end{figure}

For $0<\alpha<1,$ the harmonic mapping $f_5(z)=z+B\alpha/A_1\overline{z}+B/A_2(1-\alpha)\overline{z^2}$ is hereditarily $\lambda-$spirallike as from Lemma \ref{p3-080}, the function $F(z)=z+\alpha\overline{z}+1/2(1-\alpha)\overline{z^2}$ is hereditarily starlike. Here, we note that $d_n=nB/A_n.$ The images of $\mathbb{D}$ under $f_5$ for certain values of $\alpha$ and $\lambda$ are shown in Figure \ref{fig-3}.

Thus we see that one can construct  harmonic spirallike functions from given starlike mapping, and vice-versa. \\

As many geometric properties of analytic univalent functions are already known, it provide more geometric exposures. So, it  is quite interesting to relate harmonic function theory to analytic function theory. In this regards, Nagpal and Ravichandran \cite{Nagpal-Ravichandran-2013} established a relation between harmonic fully starlike function and analytic starlike function.
\begin{lem}\label{p3-220}\cite{Nagpal-Ravichandran-2013}
A sense-preserving harmonic function $f=h+\overline{g}$ is fully starlike in $\mathbb{D}$ if the analytic functions $h+\epsilon g$ are starlike in $\mathbb{D}$ for each $|\epsilon|=1.$
\end{lem}
Our next result for spirallike function is based on similar concept.
\begin{thm}\label{p3-230}
Let $\lambda\in(-\pi/2,\pi/2)$ and $\{d_n\}$ be sequence of complex numbers with $|d_n|\leq nB/A_n$ for $n\geq 1.$ Let
$H(z)=z-\sum\limits_{n=2}^\infty |a_n|z^n$ and  $G(z)=\sum\limits_{n=1}^\infty |b_n|z^n$
be two analytic functions. If $$F_\epsilon(z)=z\left(\frac{H(z)+\epsilon G(z)}{z}\right)^{e^{i\lambda}\cos\lambda}$$ is analytic and $\lambda-$spirallike for each $|\epsilon|=1$ then the harmonic mapping $$f(z)=h(z)+\overline{g(z)}=z+\sum\limits_{n=2}^\infty d_na_nz^n+\overline{\sum\limits_{n=1}^\infty d_nb_nz^n}$$ is hereditarily $\lambda-$spiralike.
\end{thm}

\section{Proof of Main Results}\label{Proof of Main Results}
In this section, we will proof all the results sequentially.
\begin{proof}[\textbf{Proof of Theorem \ref{p3-040}}]
Since $f\in\mathcal{SP}_H(\lambda),$ it follows that
\begin{align*}
{\rm Re\,}\left(e^{-i\lambda}\frac{Df(z)}{f(z)}\right)
&= {\rm Re\,}\left(e^{-i\lambda}\frac{zh'(z)-\overline{zg'(z)}}{h(z)+g(z)}\right)\\
& = {\rm Re\,}\left(e^{-i\lambda}\frac{\left(zh'(z)-\overline{zg'(z)}\right)(\overline{h(z)}+g(z))}{|h(z)+\overline{g(z)}|^2}\right)\\
&= \frac{1}{|h(z)+\overline{g(z)}|^2}{\rm Re\,}\left[ e^{-i\lambda}\frac{zh'(z)}{h(z)}|h(z)|^2-\overline{e^{i\lambda}\frac{zg'(z)}{g(z)}}|g(z)|^2
 \right.\\
&\hspace{40mm}\left.
   + e^{-i\lambda}zh'(z)g(z)-\overline{e^{i\lambda}zg'(z)h(z)}\right]\\
&  >0.
\end{align*}
Upon using the fact ${\rm Re\,}w={\rm Re\,}(\overline{w})$ for $w\in\mathbb{C}$ in the last inequality, we get \eqref{p3-050}. \\

Conversely, a simple calculation leads the condition \eqref{p3-050} to ${\rm Re\,}\left(e^{-i\lambda}\frac{Df(z)}{f(z)}\right)>0.$ And from the hypothesis, $f$ is sense-preserving and $f(z)=0$ only for $z=0.$ Therefore Lemma \ref{p3-030} implies that the function $f$ is harmonic hereditarily $\lambda-$spirallike function.
%If $f\in\mathcal{SP}_H(\lambda)$ then ${\rm Re\,}\left(e^{-i\lambda}\frac{Df(z)}{f(z)}\right)>0$ for $z\in\mathbb{D}$. using the fact ${\rm Re\,}w={\rm Re\,}(\overline{w})$ for $w\in\mathbb{C}$ in the last inequality, one can easily get \eqref{p3-050}. \\
%Conversely, a simple calculation leads the condition \eqref{p3-050} to ${\rm Re\,}\left(e^{-i\lambda}\frac{Df(z)}{f(z)}\right)>0.$ Furthermore, $f$ is sense-preserving and $f(z)=0$ only for $z=0.$ Therefore Lemma \ref{p3-030} implies that the function $f$ is hereditarily $\lambda-$spirallike.
\end{proof}

\begin{proof}[\textbf{Proof of Theorem \ref{p3-100}}]
We first note that
\begin{align}\label{p3-240}
\frac{A_n}{B} =\frac{\left|1+ne^{-i\lambda}\right|+\left|1-ne^{-i\lambda}\right|}{\left|1+e^{-i\lambda}\right|-\left|1-e^{-i\lambda}\right|}\geq \frac{\left|(1+ne^{-i\lambda})-(1-ne^{-i\lambda})\right|}{|1+e^{-i\lambda}|}\geq \frac{2n}{2}=n
\end{align}
and so by \eqref{p3-110} we have
\begin{align*}
|h'(z)| & \geq 1-\sum\limits_{n=2}^\infty n|a_n|r^{n-1}> 1-\sum\limits_{n=2}^\infty n|a_n|\geq 1-\sum\limits_{n=2}^\infty \frac{A_n}{B}|a_n|\\
        & \geq \sum\limits_{n=1}^\infty \frac{A_n}{B}|b_n|\geq \sum\limits_{n=1}^\infty n|b_n|>\sum\limits_{n=1}^\infty n|b_n|r^{n-1}=|g'(z)|.
\end{align*}
Thus the function $f$ is sense-preserving and locally univalent in $\mathbb{D}.$ Next, we see that
     $$\sum\limits_{n=2}^\infty n|a_n|+\sum\limits_{n=1}^\infty n|b_n|\leq \sum\limits_{n=2}^\infty \frac{A_n}{B}|a_n|+\sum\limits_{n=1}^\infty \frac{A_n}{B}|b_n|\leq1
$$
 and hence, from Lemma \ref{p3-070}, it follows that $f$ is univalent. Now, to show that $f$ is hereditarily $\lambda-$spirallike, it is sufficient to show
\begin{align*}
{\rm Re\,}\left(e^{-i\lambda}\frac{Df(z)}{f(z)}\right)>0,\quad z\in\mathbb{D}\setminus\{0\},
\end{align*}
which is equivalent to,
\begin{align*}
\left|1+e^{-i\lambda}\frac{Df(z)}{f(z)}\right|>\left|1-e^{-i\lambda}\frac{Df(z)}{f(z)}\right|,
\end{align*}
or,
\begin{align}\label{p3-250}
M:=\left|f(z)+e^{-i\lambda}Df(z)\right|-\left|f(z)-e^{-i\lambda}Df(z)\right|>0.
\end{align}
Now
\begin{align}\label{p3-260}
M &=\left|f(z)+e^{-i\lambda}Df(z)\right|-\left|f(z)-e^{-i\lambda}Df(z)\right|
\\\nonumber
&=\left|(h(z)+e^{-i\lambda}zh'(z))+\overline{(g(z)-e^{i\lambda}zg'(z))}\right|\\\nonumber
&\qquad-\left|(h(z)-e^{-i\lambda}zh'(z))+\overline{(g(z)+e^{i\lambda}zg'(z))}\right|\\\nonumber
&=\left|(1+e^{-i\lambda})z+\sum\limits_{n=2}\left(1+ne^{-i\lambda}\right)a_nz^n+\sum\limits_{n=1}\left(1-ne^{-i\lambda}\right)\overline{b_nz^n}\right|\\\nonumber
&\qquad\quad\quad-\left|(1-e^{-i\lambda})z+\sum\limits_{n=2}\left(1-ne^{-i\lambda}\right)a_nz^n+\sum\limits_{n=1}\left(1+ne^{-i\lambda}\right)\overline{b_nz^n}\right|\\\nonumber
& \geq \left(\left|1+e^{-i\lambda}\right|-\left|1-e^{-i\lambda}\right|\right)|z|-\sum\limits_{n=2}^\infty\left(\left|1+ne^{-i\lambda}\right|+\left|1-ne^{-i\lambda}\right| \right)|a_n||z|^n\\\nonumber
&\qquad\quad\quad-\sum\limits_{n=1}^\infty\left(\left|1-ne^{-i\lambda}\right|+\left|1+ne^{-i\lambda}\right| \right)|b_n||z|^n\\\nonumber
&=|z|\left[1-\left(\sum\limits_{n=2}^\infty\frac{A_n}{B}|a_n|+\sum\limits_{n=1}^\infty\frac{A_n}{B}|b_n|\right)|z|^{n-1} \right]
\end{align}
where
\begin{align*}
A_n &=\left|1+ne^{-i\lambda}\right|+\left|1-ne^{-i\lambda}\right|\quad\text{and}\quad B=\left|1+e^{-i\lambda}\right|-\left|1-e^{-i\lambda}\right|.
\end{align*}
%    &=(1+n^2+2n\cos\lambda)^{1/2}+(1+n^2-2n\cos\lambda)^{1/2}\\
%and
%\begin{align*}
%B &=\left|1+e^{-i\lambda}\right|-\left|1-e^{-i\lambda}\right|.
%\end{align*}
%  &=\sqrt{2}[(1+\cos\lambda)^{1/2}-(1-\cos\lambda)^{1/2}].
Thus from \eqref{p3-110} and \eqref{p3-260}, for $z\in\mathbb{D}\setminus\{0\},$ we have
\begin{align*}
M &\geq |z|\left[1-\left(\sum\limits_{n=2}^\infty\frac{A_n}{B}|a_n|+\sum\limits_{n=1}^\infty\frac{A_n}{B}|b_n|\right)|z|^{n-1} \right]\\
&>|z|\left[1-\left(\sum\limits_{n=2}^\infty\frac{A_n}{B}|a_n|+\sum\limits_{n=1}^\infty\frac{A_n}{B}|b_n|\right)\right]>0.
\end{align*}
\end{proof}

\begin{proof}[\textbf{Proof of Theorem \ref{p3-120}}]
First part of the theorem follows immediately from the proof the the Theorem \ref{p3-100}. So we will prove the converse part. Let $f=h+\overline{g}\in\mathcal{SP}_{H'}(\lambda).$ Then $f$ is univalent and satisfies
\begin{align}\label{p3-270}
{\rm Re\,}\left(e^{-i\lambda}\frac{Df(z)}{f(z)}\right)>0 ~~\text{for}~z\in\mathbb{D}\setminus\{0\}.
\end{align}
For $0\leq r<1,$ we define
$$\phi(r)=1-\left(\sum_{n=2}^{\infty} |a_n| -\sum_{n=1}^{\infty} |b_n|\right)r^{n-1}$$
and
$$\psi(r)=\sum_{n=2}^{\infty} n|a_n|r^{n-1} +\sum_{n=1}^{\infty} n|b_n|r^{n-1}.$$
Then on the positive real axis, $f$ can be written as $f(r)=r\phi(r),~0\leq r<1.$ Since $f(z)\neq 0$ for $z\neq 0$ it follows that $\phi(r)\neq 0$ in $(0,1).$ Indeed, we show that $\phi(r)>0$ for each $0<r<1$. On the contrary, suppose that there exist a $r_0\in (0,1)$ such that $\phi(r_0)<0$. Since $\phi$ is continuous in $[0,r_0]$ and $\phi(0)=1>0,$ it follows that $\phi$ has a zero in  $(0,r_0),$ which contradicts the fact that $f(z)\neq 0$ for $z\neq 0.$
Now,
\begin{align*}
\frac{Df(z)}{f(z)}  =\frac{zh'(z)-\overline{zg'(z)}}{h(z)+\overline{g(z)}}=\frac{z-\sum\limits_{n=2}^\infty n|a_n|z^n-\sum\limits_{n=1}^\infty n|b_n|z^n}{z-\sum\limits_{n=2}^\infty |a_n|z^n+\sum\limits_{n=1}^\infty|b_n|z^n}.
\end{align*}
Taking the values of $z$ on the positive real axis, i.e., $z=r,~0<r<1,$ we get
\begin{align*}
{\rm Re\,}\left(e^{-i\lambda}\frac{Df(r)}{f(r)}\right)= \cos\lambda\frac{1-\psi(r)}{\phi(r)}.
\end{align*}
As $f=h+\overline{g}\in\mathcal{SP}_{H'}(\lambda)$ and $\phi(r)>0,$ it follows that $1-\psi(r)>0$ for $0<r<1.$ Thus $\psi(r)<1$ in $(0,1).$ Also, we see that $\psi$ is monotonic increasing function of $r\in(0,1)$ and hence
\begin{align}\label{p3-280}
\lim_{r\rightarrow 1^{-}}\psi(r)=\sum\limits_{n=2}^\infty n|a_n|+\sum\limits_{n=1}^\infty n|b_n|)\leq1.
\end{align}
Since $B/A_n\leq n$ for $n\geq1$ (see \eqref{p3-240}), it follows that
\begin{align*}
\sum\limits_{n=2}^\infty \frac{B}{A_n}|a_n|+\sum\limits_{n=1}^\infty \frac{B}{A_n}|b_n|\leq \sum\limits_{n=2}^\infty n|a_n|+\sum\limits_{n=1}^\infty n|b_n|\leq1.
\end{align*}
This completes the proof.
%\leq \cos\lambda\frac{1-\left(\sum\limits_{n=2}^\infty \frac{B}{A_n}|a_n|r^{n-1}+\sum\limits_{n=1}^\infty\frac{B}{A_n}|b_n|r^{n-1}\right)}{1-\left(\sum\limits_{n=2}^\infty |a_n|-\sum\limits_{n=1}^\infty|b_n|\right)r^{n-1}}
%It is easy to see that the numerator is negative for $r$  very close to $1$ if \eqref{p3-130} does not hold. This completes the proof.
\end{proof}

\begin{proof}[\textbf{Proof of Corollary \ref{p3-140}}]
If $f=h+\overline{g}\in\mathcal{SP}_{H'}(\lambda)$ then by the argument given in the proof of the Theorem \ref{p3-120}, the inequality \eqref{p3-280} hold.  Now the proof follows from Lemma \ref{p3-070}.
\end{proof}

\begin{proof}[\textbf{Proof of Theorem \ref{p3-170}}]
Let $|z|=r<1.$ As the harmonic function $f(z)=z+\sum\limits_{n=2}^\infty a_nz^n+\overline{\sum\limits_{n=1}^\infty b_nz^n}$ satisfies the condition \eqref{p3-110}, it follows that
\begin{equation*}
\begin{split}
|f(z)| & \leq |z|+\sum\limits_{n=2}^\infty|a_n||z|^n+\sum\limits_{n=1}^\infty|b_n||z|^n\\
       & \leq r+\sum\limits_{n=2}^\infty|a_n|r+\sum\limits_{n=1}^\infty|b_n|r\\
       & =r+\frac{B}{A_1}\left(\sum\limits_{n=2}^\infty \frac{A_1}{B}|a_n|+\sum\limits_{n=1}^\infty \frac{A_1}{B}|b_n|\right)r\\
       & \leq r+\frac{B}{A_1}\left(\sum\limits_{n=2}^\infty \frac{A_n}{B}|a_n|+\sum\limits_{n=1}^\infty \frac{A_n}{B}|b_n|\right)r\\
       & \leq r+\frac{B}{A_1}r\\
       &= \left(1+\frac{B}{A_1}\right)r.
\end{split}
\end{equation*}
%       & \leq r+\frac{B}{A_1}\left(\sum\limits_{n=2}^\infty \frac{A_2}{B}|a_n|+\sum\limits_{n=1}^\infty \frac{A_1}{B}|b_n|\right)r\\
On the other hand,
\begin{align}\label{p3-290}
|f(z)|\nonumber & \geq |z|-\sum\limits_{n=2}^\infty|a_n||z|^n-\sum\limits_{n=1}^\infty|b_n||z|^n\\\nonumber
       & \geq r-\sum\limits_{n=2}^\infty|a_n|r-\sum\limits_{n=1}^\infty|b_n|r\\\nonumber
       & =r-\frac{B}{A_1}\left(\sum\limits_{n=2}^\infty \frac{A_1}{B}|a_n|+\sum\limits_{n=1}^\infty \frac{A_1}{B}|b_n|\right)r\\\nonumber
       & \geq r-\frac{B}{A_1}\left(\sum\limits_{n=2}^\infty \frac{A_n}{B}|a_n|+\sum\limits_{n=1}^\infty \frac{A_n}{B}|b_n|\right)r\\\nonumber
       & \geq r-\frac{B}{A_1}r\\
       &= \left(1-\frac{B}{A_1}\right)r.
\end{align}
%       & \geq r-\frac{B}{A_1}\left(\sum\limits_{n=2}^\infty \frac{A_2}{B}|a_n|+\sum\limits_{n=1}^\infty \frac{A_1}{B}|b_n|\right)r\\
Left hand estimate is sharp for the hereditarily $\lambda-$spirallike function $f_6(z)=z-B/A_1\overline{z}$ whereas the right hand estimate is sharp for $f_7(z)=z+B/A_1\overline{z}.$\\

On taking $r\rightarrow 1^{-}$ in \eqref{p3-290}, we get
%$f(z)\geq \left(1-\frac{B}{A_1}\right).$ Thus
\begin{align*}
\left\{w\in\mathbb{C}:|w|<1-\frac{B}{A_1}\right\}\subset f(\mathbb{D}).
\end{align*}

\end{proof}

\begin{proof}[\textbf{Proof of Theorem \ref{p3-200}}]
As $h(z)/z$ and $g(z)/z$ are non zero in the unit disk $\mathbb{D},$ from the relation \eqref{p3-190}, we have
  $$e^{i\lambda}\log\frac{h(z)}{z}=\cos\lambda\log\frac{g(z)}{z}.$$
Upon taking the differential operator D, we get
$$e^{i\lambda}\left(\frac{Dh(z)}{h(z)}-1\right)=\cos\lambda\left(\frac{Dg(z)}{g(z)}-1\right)$$
which is equivalent to
$$e^{i\lambda}\frac{Dh(z)}{h(z)}=\cos\lambda\frac{Dg(z)}{g(z)}+i\sin\lambda.$$
Thus $${\rm Re\,}\left(e^{i\lambda}\frac{Dh(z)}{h(z)}\right)=\cos\lambda{\rm Re\,}\left(\frac{Dg(z)}{g(z)}\right)$$
from which the result follows.
\end{proof}

\begin{proof}[\textbf{Proof of Theorem \ref{p3-210}}]
Since $F(z)=z-\sum\limits_{n=2}^\infty |a_n|z^n+\overline{\sum\limits_{n=1}^\infty |b_n|z^n}$ is harmonic hereditarily starlike mapping of the form \eqref{p3-020}, by Lemma \ref{p3-070}, the coefficient condition \eqref{p3-080} holds.
%\begin{align}\label{pp-320}
%\sum\limits_{n=1}^\infty n(|a_n|+|b_n|)\leq 2.
%\end{align}
Since $\{d_n\}$ is a sequence of complex numbers with $|d_n|\leq nB/A_n$ for all $n\geq 1,$ then
\begin{align*}
\sum\limits_{n=2}^\infty\frac{A_n}{B}|d_na_n|+\sum\limits_{n=1}^\infty\frac{A_n}{B}|d_nb_n| & = \sum\limits_{n=2}^\infty\frac{A_n}{B}|d_n||a_n|+\sum\limits_{n=1}^\infty\frac{A_n}{B}|d_n||b_n|\\
  & \leq \sum\limits_{n=2}^\infty n|a_n|+\sum\limits_{n=1}^\infty n|b_n|\\
  & \leq 1.
\end{align*}
Hence by Theorem \ref{p3-100}, $f(z)=z+\sum\limits_{n=2}^\infty d_n a_nz^n+\overline{\sum\limits_{n=1}^\infty d_n b_nz^n}$ is hereditarily $\lambda-$spirallike.
Conversely, if $f(z)=z-\sum\limits_{n=2}^\infty |a_n|z^n+\overline{\sum\limits_{n=1}^\infty |b_n|z^n}$ is harmonic hereditarily $\lambda-$spirallike mapping, then from the Corollary \ref{p3-140}, the coefficient condition \eqref{p3-080} holds. Hence by Lemma \ref{p3-070}, the function $F(z)=z+\sum\limits_{n=2}^\infty a_n z^n+\overline{\sum\limits_{n=1}^\infty b_n z^n}$ is hereditarily starlike.
\end{proof}

\begin{proof}[\textbf{Proof of Theorem \ref{p3-230}}]
As $F_\epsilon(z)=z\left(\frac{H(z)+\epsilon G(z)}{z}\right)^{e^{i\lambda}\cos\lambda}$ is $\lambda-$spirallike  for each $|\epsilon|=1$, from the Lemma \ref{p3-180}, it follows that $H(z)+\epsilon G(z)$ is analytic starlike for each $|\epsilon|=1.$ Consequently, by Lemma \ref{p3-220}, the associated harmonic function $F(z)=H(z)+\overline{G(z)}$ is fully starlike  of the form \eqref{p3-020}. Thus, from Theorem \ref{p3-210}, $f=z+\sum\limits_{n=2}^\infty d_na_nz^n+\overline{\sum\limits_{n=1}^\infty d_nb_nz^n}$ is harmonic hereditarily $\lambda-$spirallike.
%Thus $f(z)=z+\sum\limits_{n=2}^\infty d_na_nz^n+\overline{\sum\limits_{n=1}^\infty d_nb_nz^n}\in\mathcal{SP}_H(\lambda)$ by Theorem \ref{p3-100}.
\end{proof}

\noindent\textbf{Declarations:\\}

\noindent\textbf{Data availability:}
Data sharing not applicable to this article as no data sets were generated or analyzed during the current study.\\

\noindent\textbf{Authors Contributions:}
All authors contributed equally to the investigation of the problem and the order of the authors is given alphabetically according to their surname. All authors read and approved the final manuscript. \\

\noindent\textbf{Acknowledgement:}
The second named author thanks the Department Of Science and Technology, Ministry Of Science and Technology, Government Of India
for the financial support through DST-INSPIRE Fellowship (No. DST/INSPIRE Fellowship/2018/IF180967).\\

\textbf{Conflict of interest:} The authors declare that they have no conflict of interest.

\end{document}